# THE SHAPE OF INCOMPLETE PREFERENCES[1]

By Robert Nau

*Duke University*

Incomplete preferences provide the epistemic foundation for models of imprecise subjective probabilities and utilities that are used in robust Bayesian analysis and in theories of bounded rationality. This paper presents a simple axiomatization of incomplete preferences and characterizes the shape of their representing sets of probabilities and utilities. Deletion of the completeness assumption from the axiom system of Anscombe and Aumann yields preferences represented by a convex set of state-dependent expected utilities, of which at least one must be a probability/utility pair. A strengthening of the state-independence axiom is needed to obtain a representation purely in terms of a set of probability/utility pairs.

**1. Introduction.** In the Bayesian theory of choice under uncertainty, a decision maker holds rational preferences among acts that are mappings from states of nature $\{s\}$ to consequences $\{c\}$. It is typically assumed that rational preferences are *complete*, meaning that, for any two acts $\mathbf{X}$ and $\mathbf{Y}$, either $\mathbf{X} \succsim \mathbf{Y}$ ($\mathbf{X}$ is weakly preferred to $\mathbf{Y}$) or else $\mathbf{Y} \succsim \mathbf{X}$, or both. This assumption, together with other rationality axioms such as transitivity and independence, leads to a representation of preferences by a unique subjective probability distribution on states $p(s)$ and a unique utility function $u(c)$ on consequences, such that $\mathbf{X} \succsim \mathbf{Y}$ if and only if the subjective expected utility of $\mathbf{X}$ is greater than or equal to that of $\mathbf{Y}$ [1, 5, 24]. However, the completeness assumption may be inappropriate if only partial information about the decision maker's preferences is available, or if realistic limits on her powers of discrimination are assumed, or if there are many decision makers whose preferences may disagree.

Received September 2003; revised October 2005.
[1]Supported in part by the Fuqua School of Business and by the National Science Foundation.

Incomplete preferences are generally represented by imprecise (set-valued) probabilities and/or utilities. Varying degrees and types of imprecision have been modeled previously in the literature of statistical decision theory and rational choice:

1. If probabilities alone are considered to be imprecise, then preferences can be represented by a set of probability distributions $\{p(s)\}$ and together with a unique, exogenously-specified utility function $u(c)$. The set of probability distributions is typically convex, so the representation can be derived by separating hyperplane arguments (e.g., [7, 15, 29, 31, 33]). Representations of this kind are widely used in robust Bayesian statistics [22, 32] and they are also attracting interest in economics [19].
2. If utilities alone are considered to be imprecise, preferences can be represented by a set of utility functions $\{u(c)\}$ and a unique, exogenously-specified probability distribution $p(s)$, a representation that has been axiomatized and applied to economic models by Aumann [2] and Dubra, Maccheroni and Ok [4]. The set of utility functions in this case is also typically convex, so that separating hyperplane arguments are again applicable: the roles of probabilities and utilities are merely reversed.
3. If both probabilities and utilities are allowed to be imprecise, they can be represented by separate sets of probability distributions $\{p(s)\}$ and utility functions $\{u(c)\}$ whose elements are paired up arbitrarily. This representation of preferences preserves the traditional separation of information about *beliefs* from information about *values* when both are imprecise [20, 21]. It lacks a compelling axiomatic basis, but it arises naturally when imprecise probabilities and utilities are assessed independently, as they often are in practice.
4. More generally, incomplete preferences can be represented by sets of probability distributions paired with state-independent utility functions $\{(p(s), u(c))\}$, a.k.a. "probability/utility pairs." This representation has an appealing multi-Bayesian interpretation and provides a normative basis for techniques of robust decision analysis [14] and asset pricing in incomplete financial markets [30]. It has been axiomatized by Seidenfeld, Schervish and Kadane [28] (henceforth SSK), starting from the "horse lottery" formalization of decision theory introduced by Anscombe and Aumann [1]. However, as pointed out by SSK, the set of probability/utility pairs is typically nonconvex and may even be unconnected, so that separating hyperplane arguments are hard to apply. Instead, SSK introduce methods of transfinite induction and indirect reasoning.

The objective of this paper is to derive a simple representation of incomplete preferences for the elementary case of finite state and reward spaces, and to characterize the shape of the resulting sets of probabilities and utilities. First, it is shown that deleting both completeness and state-independence from the horse-lottery axiom system of Anscombe and Aumann leads to a representation of preferences by a set of probabilities paired with state-*dependent* utility functions $\{(p(s), u(s,c))\}$. Such pairs will be called *state-dependent expected utility* (s.d.e.u.) functions. State-dependent utilities have been used in economic models by Karni [9] and Drèze [3] and are also discussed by Schervish, Seidenfeld and Kadane [25]. A set of s.d.e.u. functions is typically convex—unlike a set of probability/utility pairs—so that separating-hyperplane methods are still applicable at this stage. Next, the usual state-independence axiom is reintroduced and shown to impose (only) the further requirement that the representing set should contain *at least one* probability/utility pair, analogous to SSK's result establishing one-way agreement between a partially ordered preference relation and a nonempty set of probability/utility pairs. However, those assumptions are too weak to yield two-way agreement in which every extremal preference has an agreeing state-independent utility, because the state-independence axiom may have limited applicability when preferences are incomplete. The main result of this paper is to show that two-way agreement between incomplete preferences and state-independent utilities can be achieved by strengthening the state-independence axiom. Although the representing set of probability/utility pairs is nonconvex, it is characterized simply as the intersection of a convex set of s.d.e.u. functions and the nonconvex set of state-independent utilities.

The organization of the paper is as follows. Section 2 introduces basic notation and derives a representation of preferences by convex sets of s.d.e.u. functions when neither completeness nor state-independence is assumed. Section 3 incorporates Anscombe and Aumann's state-independence assumption and shows that it requires (only) the existence of at least one agreeing state-independent utility. Section 4 discusses an example of SSK to highlight the implications of different continuity and strictness conditions. Section 5 presents the stronger constructive axiom that is needed to obtain a representation purely in terms of probability/utility pairs, illustrated by another example. Section 6 briefly discusses the results. Proofs are given in the Appendix.

**2. Representation of incomplete preferences.** Let $S$ denote a finite set of states and let $C$ denote a finite set of consequences. Let $\mathcal{B} = \{\mathbf{B} : S \times C \mapsto \Re\}$. An element $\mathbf{X} \in \mathcal{B}$ is a (*horse*) *lottery* if $\mathbf{X} \geq \mathbf{0}$ and $\forall s, \sum_c X(s,c) = 1$, with the interpretation that $X(s,c)$ is the objective probability of receiving consequence $c$ when state $s$ occurs. Henceforth, the symbols $\mathbf{X}$, $\mathbf{Y}$, $\mathbf{Z}$ and $\mathbf{H}$ will be used to denote lotteries; the symbol $\mathbf{B}$ will denote an element of $\mathcal{B}$ that is not necessarily a lottery (e.g., $\mathbf{B}$ will typically represent the difference between two lotteries). A lottery $\mathbf{X}$ is *constant* if the probabilities it assigns to consequences are constant across states—that is, if $X(s,c) = X(s',c)$ for all $s, s', c$. The symbol $\succsim$ will denote weak preference between lotteries: $\mathbf{X} \succsim \mathbf{Y}$ means that $\mathbf{X}$ is preferred or indifferent to $\mathbf{Y}$, which is the behavioral primitive. The domain of $\succsim$ is the set of all lotteries. The asymmetric part of $\succsim$ is denoted by $\succ$, but it will not be used except in stating axiom A5. An *extension* of $\succsim$ is any preference relation $\succsim_0$ that is stronger in the sense that $\mathbf{X} \succsim \mathbf{Y}$ implies $\mathbf{X} \succsim_0 \mathbf{Y}$ but not necessarily vice versa.

An *event* is a subset of $S$. The symbols $\mathbf{E}$ and $\mathbf{F}$ will be used interchangeably as names for events and for their indicator functions on $S \times C$, so that for all $c$, $E(s,c) = 1[0]$ if the event $\mathbf{E}$ includes [does not include] state $s$. $\mathbf{E}_s$ will denote the indicator vector for state $s$, so that for all $c$, $\mathbf{E}_s(s',c) = 1[0]$ if $s = s'[s \neq s']$. If $\alpha$ is a scalar between 0 and 1, then $\alpha \mathbf{X} + (1-\alpha)\mathbf{Y}$ is an *objective mixture* of $\mathbf{X}$ and $\mathbf{Y}$: it yields consequence $c$ in state $s$ with probability $\alpha X(s,c) + (1-\alpha)Y(s,c)$. If $\mathbf{E}$ is an event, then $\mathbf{EX} + (\mathbf{1} - \mathbf{E})\mathbf{Y}$ is a *subjective mixture* of $\mathbf{X}$ and $\mathbf{Y}$: it yields consequence $c$ in state $s$ with probability $X(s,c)$ if $E(s,c) = 1$, and with probability $Y(s,c)$ otherwise.

Assume that $C$ contains a "worst" and a "best" consequence, labeled 0 and 1, respectively. This assumption follows Luce and Raiffa [13] and Anscombe and Aumann [1], and it is technically without loss of generality in the sense that the preference order can always be extended to a larger domain that includes two additional consequences which, by construction, are better and worse, respectively, than all the original consequences. (Such an extension is demonstrated by SSK, Theorem 2.) The best and worst consequences ultimately serve to calibrate the definition and measurement of subjective probabilities, but even so the probabilities remain somewhat arbitrary, as will be shown. Intermediate consequences are labeled $2, 3, \ldots, K$, and for all $c \in \{0, 1, 2, \ldots, K\}$, $\mathbf{H}_c$ denotes the special lottery that yields consequence $c$ with probability 1 in every state. That is, for all $s$, $\mathbf{H}_c(s,c') = 1[0]$ if $c' = c[c' \neq c]$.

The first group of axioms that are assumed to govern rational preference is as follows:

A1 (*Quasi order*) $\succsim$ is transitive and reflexive.
A2 (*Mixture-independence*) $\mathbf{X} \succsim \mathbf{Y} \Leftrightarrow \alpha\mathbf{X} + (1-\alpha)\mathbf{Z} \succsim \alpha\mathbf{Y} + (1-\alpha)\mathbf{Z} \,\forall\, \alpha \in (0,1)$.
A3 (*Continuity in probability*) If $\{\mathbf{X}_n\}$ and $\{\mathbf{Y}_n\}$ are convergent sequences such that $\mathbf{X}_n \succsim \mathbf{Y}_n$, then $\lim \mathbf{X}_n \succsim \lim \mathbf{Y}_n$.
A4 (*Existence of best and worst*) For all $c > 1$, $\mathbf{H}_1 \succsim \mathbf{H}_c \succsim \mathbf{H}_0$.
A5 (*Coherence or nontriviality*) $\mathbf{H}_1 \succ \mathbf{H}_0$ (i.e., *not* $\mathbf{H}_0 \succsim \mathbf{H}_1$).

A1 and A2 are von Neumann and Morgenstern's first two axioms of expected utility, minus completeness, as applied to horse lotteries by Anscombe and Aumann [1]; see also [5]. A3 is a strong continuity condition used by Garcia del Amo and Ríos Insua [6] that also works in infinite-dimensional spaces.

A4 and A5 ensure nontriviality and provide reference points for probability measurement, as noted earlier.

The following definition (following SSK) will prove useful for distinguishing a subset of preferences that implicitly determines the entire preference order:

DEFINITION. A collection of preferences $\{\mathbf{X}_n \succsim \mathbf{Y}_n\}$ is a *basis* for $\succsim$ under an axiom system if every preference $\mathbf{X} \succsim \mathbf{Y}$ can be deduced from $\{\mathbf{X}_n \succsim \mathbf{Y}_n\}$ by direct application of those axioms.

The primal geometric representation of $\succsim$ is now given by the following:

THEOREM 1. $\succsim$ *satisfies* A1–A5 *if and only if there exists a closed convex cone* $\mathcal{B}^* \subset \mathcal{B}$, *receding from the origin, such that* $\mathbf{X} \succsim \mathbf{Y} \Leftrightarrow \mathbf{X} - \mathbf{Y} \in \mathcal{B}^*$. *In particular, if* $\{\mathbf{X}_n \succsim \mathbf{Y}_n\}$ *is a basis for* $\succsim$ *under* A1–A5, *then* $\mathcal{B}^*$ *is the closed convex hull of the rays whose directions are* $\{\mathbf{X}_n - \mathbf{Y}_n\}$ *for all* $n$, *together with* $\{\mathbf{H}_1 - \mathbf{H}_c\}$ *and* $\{\mathbf{H}_c - \mathbf{H}_0\}$ *for all* $c$.

Thus, the direction of preference between two lotteries $\mathbf{X}$ and $\mathbf{Y}$ depends only on their difference $\mathbf{X} - \mathbf{Y}$, with $\mathbf{X} \succsim \mathbf{Y}$ if and only if $\mathbf{X} - \mathbf{Y}$ is in a convex cone $\mathcal{B}^*$ of *preferred directions*. It follows that if $\mathbf{EX} + (\mathbf{1} - \mathbf{E})\mathbf{Z} \succsim \mathbf{EY} + (\mathbf{1} - \mathbf{E})\mathbf{Z}$, where $\mathbf{E}$ is an event, then $\mathbf{EX} + (\mathbf{1} - \mathbf{E})\mathbf{Z}' \succsim \mathbf{EY} + (\mathbf{1} - \mathbf{E})\mathbf{Z}'$ for any $\mathbf{Z}'$, which is the so-called "sure-thing principle": where two lotteries agree, it does not matter *how* they agree there. Henceforth, because agreeing conditional components of lotteries do not affect preferences, the expression $\mathbf{EX} \succsim \mathbf{EY}$ will be used as shorthand for $\mathbf{EX} + (\mathbf{1} - \mathbf{E})\mathbf{Z} \succsim \mathbf{EY} + (\mathbf{1} - \mathbf{E})\mathbf{Z}$ for all $\mathbf{Z}$, which means that $\mathbf{X}$ is preferred to $\mathbf{Y}$ *conditional on the event* $\mathbf{E}$, or, equivalently, that $\mathbf{E}(\mathbf{X} - \mathbf{Y})$ is a preferred direction. Similarly, for conditional lotteries embedded in objective mixtures, $\alpha \mathbf{EX} + (1-\alpha)\mathbf{X}' \succsim \alpha \mathbf{EY} + (1-\alpha)\mathbf{Y}'$ will be used (in Section 5) as shorthand for $\alpha(\mathbf{EX} + (\mathbf{1} - \mathbf{E})\mathbf{Z}) + (1-\alpha)\mathbf{X}' \succsim \alpha(\mathbf{EY} + (\mathbf{1} - \mathbf{E})\mathbf{Z}) + (1-\alpha)\mathbf{Y}'$ for all $\mathbf{Z}$, meaning that $\alpha \mathbf{E}(\mathbf{X} - \mathbf{Y}) + (1-\alpha)(\mathbf{X}' - \mathbf{Y}')$ is a preferred direction.

To set the stage for the representation of incomplete preferences by sets of state-dependent utilities, let a *state-dependent expected utility* (*s.d.e.u.*) *function* be defined as a function $v: S \times C \mapsto \Re$, with the interpretation that $v(s, c)$ is an expected utility associated with consequence $c$ when it is to be received in state $s$, and let $U_v: \mathcal{B} \mapsto \Re$ denote the utility function on lotteries that is induced by $v$ according to the linear formula

$$U_v(\mathbf{X}) \equiv \sum_{s \in S, c \in C} X(s, c) v(s, c).$$

In particular, $U_v(\mathbf{E}_s \mathbf{H}_c + (\mathbf{1} - \mathbf{E}_s)\mathbf{H}_0) = v(s, c)$.

DEFINITIONS. (i) An s.d.e.u. function $v$ is a *probability/utility pair* if it can be expressed as the product of a probability distribution $p\colon S \mapsto [0,1]$ and a state-independent utility function $u\colon C \mapsto \Re$, so that $v(s,c) = p(s)u(c)$ for all $s$ and $c$. (ii) An s.d.e.u. function $v$ *agrees* (one way) with $\succsim$ if $\mathbf{X} \succsim \mathbf{Y} \Rightarrow U_v(\mathbf{X}) \geq U_v(\mathbf{Y})$. (iii) A set $\mathcal{V}$ of s.d.e.u. functions *represents* $\succsim$ if $\mathbf{X} \succsim \mathbf{Y} \Leftrightarrow U_v(\mathbf{X}) - U_v(\mathbf{Y}) \geq 0 \ \forall v \in \mathcal{V}$.

A $v$ that agrees with $\succsim$ is unique only up to positive linear scaling and the addition of state-dependent constants and, by A4, it must satisfy $U_v(\mathbf{H}_0) \leq U_v(\mathbf{H}_c) \leq U_v(\mathbf{H}_1)$; hence, there is no loss of generality in assuming that it belongs to the normalized set

$$\mathcal{V}^+ \equiv \left\{ v : v(s,0) = 0 \ \forall s;\ 0 \leq \sum_{s \in S} v(s,c) \leq 1 \ \forall c \geq 2;\ \sum_{s \in S} v(s,1) = 1 \right\}.$$

In these terms, the dual to Theorem 1 can now be given as follows:

THEOREM 2. $\succsim$ *satisfies* A1–A5 *if and only if it is represented by a nonempty closed convex set of s.d.e.u. functions* $\mathcal{V}^* \in \mathcal{V}^+$. *In particular, if* $\{\mathbf{X}_n \succsim \mathbf{Y}_n\}$ *is a basis for* $\succsim$, *then* $\mathcal{V}^*$ *is the set of* $v \in \mathcal{V}^+$ *satisfying* $\{U_v(\mathbf{X}_n) \geq U_v(\mathbf{Y}_n)\}$.

The proof relies on a separating hyperplane argument. (For a similar result on a more general space, see [21].) If the basis is finite, then $\mathcal{V}^*$ is a convex polytope whose elements need not be probability/utility pairs. Subsequent sections of the paper will discuss the additional assumptions needed to ensure that some points of $\mathcal{V}^*$—especially its extreme points—are probability/utility pairs.

**3. The state-independence axiom.** An additional axiom of Anscombe and Aumann, which they call "monotonicity in the prizes" or "substitutability," provides the usual separation of subjective probability from utility, and its implications in the context of incompleteness will now be explored. (Essentially the same axiom appears as P3 in Savage's system [24].) Some additional notation will be helpful. First, for all $p \in (0,1)$, let $\mathbf{H}_p$ denote the lottery that yields the best and worst consequences with probabilities $p$ and $1-p$ in every state, which is the objective mixture

$$\mathbf{H}_p \equiv p\mathbf{H}_1 + (1-p)\mathbf{H}_0.$$

Next define $\mathbf{H_E}$ as the lottery that yields the best consequence if event $\mathbf{E}$ occurs and the worst consequence otherwise, that is, the subjective mixture

$$\mathbf{H_E} \equiv \mathbf{E}\mathbf{H}_1 + (\mathbf{1} - \mathbf{E})\mathbf{H}_0.$$

A not-potentially-null event can then be defined by a comparison of subjective and objective mixtures:

DEFINITION. An event $\mathbf{E}$ is *not potentially null* if $\mathbf{H_E} \succsim \mathbf{H}_p$ for some $p > 0$.

The additional axiom asserts that conditional preferences among *constant* lotteries may be propagated from not-potentially-null events to all other conditioning events:

A6 (*State-independence*) If $\mathbf{X}$ and $\mathbf{Y}$ are constant and $\mathbf{E}$ is not potentially null, then $\mathbf{EX} \succsim \mathbf{EY} \Rightarrow \mathbf{FX} \succsim \mathbf{FY}$ for every other event $\mathbf{F}$.

This assumption of state-independent *preferences* among constant lotteries makes possible the assignment of state-independent *utilities* to the underlying consequences, which is central to the program of defining subjective probabilities in terms of preferences. A6 together with A4 ensures that consequences 0 and 1 are worst and best in every state, that is, $\mathbf{E}_s\mathbf{H}_0 \precsim \mathbf{E}_s\mathbf{H_c} \precsim \mathbf{E}_s\mathbf{H}_1$ for all $s$ and $c$, which implies $U_v(\mathbf{E}_s\mathbf{H}_0) \leq U_v(\mathbf{E}_s\mathbf{H_c}) \leq U_v(\mathbf{E}_s\mathbf{H}_1)$ for every agreeing $v$, which is equivalent to $v(s, 0) \leq v(s, c) \leq v(s, 1)$. The s.d.e.u. functions representing a relation that satisfies A6 can therefore be normalized to lie within the set

$$\mathcal{V}^{++} \equiv \left\{ v : 0 = v(s, 0) \leq v(s, c) \leq v(s, 1) \leq 1 \right.$$

$$\left. \forall s \in S, c \geq 2; \sum_{s \in S} v(s, 1) = 1 \right\}.$$

Under this normalization, it is "as if" consequences 0 and 1 have state-independent utilities of 0 and 1, respectively, and the *probability* assigned to state $s$ by $v$ can be interpreted as

$$p_v(s) \equiv p_v(\mathbf{E}_s) = v(s, 1),$$

because this is the expected utility of a lottery that yields a utility of 1 if state $s$ obtains and a utility of 0 otherwise. Correspondingly, the probability assigned to any event $\mathbf{E}$ by $v$ is defined by the summation

$$p_v(\mathbf{E}) \equiv U_v(\mathbf{H_E}) = \sum_{s \in \mathbf{E}} p_v(s).$$

(The same approach is used by Karni [10].) Under this interpretation, bounds on subjective probabilities can be expressed by the decision maker in either primal or dual terms in light of Theorems 1–2. The assertion that "the probability of $\mathbf{E}$ is at least $p$," that is, that $p$ is a *lower probability* for $\mathbf{E}$, is defined primally by the preference $\mathbf{H_E} \succsim \mathbf{H}_p$ and dually by the constraint that $p_v(\mathbf{E}) \geq p$ for any $v$ agreeing with $\succsim$. Upper probabilities are defined analogously. An event is not potentially null if it has a strictly positive lower probability.

Unfortunately, the attribution of state-independent utilities to consequences is ultimately arbitrary and untestable: the decision maker's true utilities could have state-dependent origin and scale factors, even if the state-independence axiom is satisfied, in which case the conventionally-defined subjective probabilities would not represent true degrees of belief. (State-dependent utility scale factors would have exactly the same effects as beliefs.) The classic behavioral definitions of subjective probability given by Savage, Anscombe–Aumann and others all suffer from the same arbitrariness. The intrinsic impossibility of inferring true probabilities from material preferences is discussed in more depth by Karni, Schmeidler and Vind [12], Rubin [23], Kadane and Winkler [8], Schervish, Seidenfeld and Kadane [25], Karni and Mongin [11] and Nau [16, 17].

Notwithstanding those caveats, the s.d.e.u. function $v$ can be further decomposed by dividing the expected utility of consequence $c$ in state $s$ by the probability of the state to obtain

$$u_v(s,c) \equiv v(s,c)/p_v(s) \qquad \text{if } p_v(s) > 0$$

as the utility assigned to consequence $c$ in state $s$. This utility is state-independent if $v$ is a probability/utility pair, and otherwise it is state-dependent. In these terms, the expected utility assigned to any lottery $\mathbf{X}$ by $v$ can be rewritten as the expected value of a possibly-state-dependent utility:

$$U_v(\mathbf{X}) = \sum_s p_v(s) \sum_c u_v(s,c) X(s,c).$$

If $v$ is a probability/utility pair, then $v(s,c) = p_v(\mathbf{E}_s) U_v(\mathbf{H}_c)$ for all $s$ and $c$. Bounds on expected utilities can also be expressed in both primal and dual terms. The assertion that "the expected utility of $\mathbf{X}$ is at least $u$" is defined primally by the preference $\mathbf{X} \succsim \mathbf{H}_u$ and dually by the constraint that $U_v(\mathbf{X}) \geq u$ for any $v$ agreeing with $\succsim$. Similarly, for conditional expected utility, define

$$U_v(\mathbf{X}|\mathbf{E}) \equiv U_v(\mathbf{X}\mathbf{E})/p_v(\mathbf{E}),$$

provided $p_v(\mathbf{E}) > 0$. Then the assertion "the conditional expected utility of $\mathbf{X}$ given $\mathbf{E}$ is at least $u$" has the primal definition $\mathbf{E}\mathbf{X} \succsim \mathbf{E}\mathbf{H}_u$ and the dual definition that $U_v(\mathbf{X}|\mathbf{E}) \geq u$ for any $v$ agreeing with $\succsim$ and satisfying $p_v(\mathbf{E}) > 0$, because, for any agreeing $v$

$$\mathbf{E}\mathbf{X} \succsim \mathbf{E}\mathbf{H}_u \implies U_v(\mathbf{E}\mathbf{X}) \geq U_v(\mathbf{E}\mathbf{H}_u) = u p_v(\mathbf{E})$$
$$\iff U_v(\mathbf{X}|\mathbf{E}) \geq u \text{ or else } p_v(\mathbf{E}) = 0.$$

If preferences are complete, in addition to A1–A5, then the primal representation $\mathcal{B}^*$ expands to an entire half-space of preferred directions and the dual representation $\mathcal{V}^*$ shrinks to a unique s.d.e.u. function $v^*$. When A6 is also assumed, the unique $v^*$ must be a probability/utility pair, which is the result obtained by Anscombe and Aumann [1]. In the absence of completeness, the contribution of A6 to the separation of probability and utility is weaker, as summarized by the main theorem of this section:

THEOREM 3. $\succsim$ *satisfies* A1–A6 *only if it is represented by a nonempty closed convex set of s.d.e.u. functions* $\mathcal{V}^{**} \in \mathcal{V}^{++}$ *on which the maximum and minimum expected utilities of every constant lottery are achieved at probability/utility pairs. In particular, if* $\{\mathbf{X}_n \succsim \mathbf{Y}_n\}$ *is a basis for* $\succsim$ *under axioms* A1–A6, *then* $\mathcal{V}^{**}$ *contains every probability/utility pair* $v \in \mathcal{V}^{++}$ *that satisfies* $\{U_v(\mathbf{X}_n) \geq U_v(\mathbf{Y}_n)\}$, *of which there must be at least one.*

Apart from the fact that $\mathcal{V}^{**}$ contains all the probability/utility pairs in $\mathcal{V}^{++}$ that agree with the basis preferences, its "shape" is not easy to describe, as will be illustrated in Section 5.

An immediate implication of Theorem 3 is the property of *stochastic dominance*, namely, that if $\mathbf{X}$ is obtained from $\mathbf{Y}$ by shifting probability mass to consequence 1 from any other consequence, and/or from consequence 0 to any other consequence, in any state, then $\mathbf{X} \succsim \mathbf{Y}$ because in this case $U_v(\mathbf{X}) \geq U_v(\mathbf{Y})$ for all $v \in \mathcal{V}^{++}$. To make this result more precise, let the $[\cdot]_{\min}$ (minimum s.d.e.u.) operation be defined on $\mathcal{B}$ as

$$[\mathbf{B}]_{\min} \equiv \min_{v \in \mathcal{V}^{++}} U_v(\mathbf{B}) = \min_{s \in S}\left[B(s,1) + \sum_{c \geq 2} \min\{0, B(s,c)\}\right].$$

This quantity is the minimum possible state-dependent expected utility that could be assigned to $\mathbf{B}$: it is achieved by assigning, within each state, a utility of 0 to those consequences $c \geq 2$ for which $\mathbf{B}$ is positive and a utility of 1 to those consequences $c \geq 2$ for which $\mathbf{B}$ is negative, then assigning a subjective probability of 1 to the state in which the conditional expected utility of $\mathbf{B}$ is minimized. Stochastic dominance and the negative orthant in $\mathcal{B}$ can now be defined in a natural way:

DEFINITIONS. (i) $\mathbf{X} \geq^* [>^*]\mathbf{Y}$ ("$\mathbf{X}$ [*strictly*] *dominates* $\mathbf{Y}$") if $[\mathbf{X} - \mathbf{Y}]_{\min} \geq [>]0$. (ii) The *open negative orthant* $\mathcal{B}^-$ consists of those $\mathbf{B}$ that are strictly dominated by the zero vector, that is, $\mathcal{B}^- = \{\mathbf{B} \in \mathcal{B} : \mathbf{0} >^* \mathbf{B}\} = \{\mathbf{B} \in \mathcal{B} : B(s,1) + \sum_{c \geq 2} \max\{0, B(s,c)\} < 0 \,\forall\, s\}$.

Theorem 3 then implies that $\mathbf{X} \geq^* [>^*]\mathbf{Y} \Rightarrow \mathbf{X} \succsim [\succ]\mathbf{Y}$.

**4. Strict vs. weak preference: an example.** The closedness of the representative sets of s.d.e.u. functions in Theorems 2 and 3 is attributable to the use of weak preference as the behavioral primitive, together with a strong continuity assumption. In contrast, SSK use strict preference as the behavioral primitive, together with a weaker continuity assumption, to explicitly allow for the representation of incomplete preferences by open sets that may fail to contain probability/utility pairs.

The differences in these approaches are illustrated by an example of SSK (Example 4.1) involving two states and three consequences, that is, $S = \{1, 2\}$ and $C = \{0, 1, 2\}$. Consequences 0 and 1 have state-independent utilities of 0 and 1 by assumption, so that a probability/utility pair can be parameterized by the probability assigned to state 1 and the utility assigned to consequence 2. Consider two probability/utility pairs $(p_i, u_i)$ in which $p_0(1) = 0.1$, $p_1(1) = 0.3$, $u_0(2) = 0.1$ and $u_1(2) = 0.4$. Let $v_0$ and $v_1$ denote the corresponding s.d.e.u. functions—that is, $v_i(s, c) = p_i(s)u_i(c)$ for $i = 0, 1$. Then $U_{v_i}(\mathbf{X})$ denotes the expected utility assigned to lottery $\mathbf{X}$ by $(p_i, u_i)$, with $U_{v_0}(\mathbf{H}_2) = 0.1$ and $U_{v_1}(\mathbf{H}_2) = 0.4$. Now let $\succ$ be defined as the preference relation that satisfies a weak Pareto condition with respect to these two probability/utility pairs—that is, $\mathbf{X} \succ \mathbf{Y} \Leftrightarrow \{U_{v_0}(\mathbf{X}) > U_{v_0}(\mathbf{Y})$ and $U_{v_1}(\mathbf{X}) > U_{v_1}(\mathbf{Y})\}$. Any s.d.e.u. function that is a convex combination of $v_0$ and $v_1$ also agrees with $\succ$, so the representing set $\mathcal{V}^{**}$ is the closed line segment whose endpoints are $v_0$ and $v_1$, but none of its interior points are probability/utility pairs.

Next SSK extend $\succ$ to obtain a new preference relation $\succ''$ by imposing the additional strict preferences $\mathbf{H}_{0.4} \succ'' \mathbf{H}_2 \succ'' \mathbf{H}_{0.1}$. The effect of this extension is to remove the two endpoints of the representing set of s.d.e.u. functions, so that $\succ''$ is represented by the *open* line segment connecting $v_0$ with $v_1$. SSK point out that, although $\succ''$ satisfies all their axioms, there is no agreeing probability/utility pair for it, since the only two candidates have been deliberately excluded. They proceed to axiomatize the concept of "almost state-independent" utilities, which agree with a strict preference relation and are "within $\varepsilon$" of being state-independent. Clearly, $\succ''$ has an almost-state-independent representation, containing points arbitrarily close to $v_0$ and $v_1$.

In the present framework, where weak preference is primitive, there is no way to implement a constraint such as $\mathbf{H}_2 \succ \mathbf{H}_{0.1}$ except by asserting that $\mathbf{H}_2 \succsim \mathbf{H}_{0.1+\varepsilon}$ for some $\varepsilon > 0$. If this assertion is made, axiom A6 begins to nibble on the $v_0$ end of the line segment and continues nibbling until the representation collapses to the $v_1$ end. (See [18] for details.) If instead the other endpoint is removed, by adding the constraint $\mathbf{H}_{0.4-\varepsilon} \succsim \mathbf{H}_2$ for $\varepsilon > 0$, collapse occurs to $v_0$. If both constraints are added, the entire interval is annihilated, yielding incoherence (a violation of A5). Hence, this example is unstable in the sense that any *finite* extension of the original preference relation leads to a collapse to one or the other of the original probability/utility pairs, or else to incoherence.

**5. The need for stronger state-independence.** In the preceding example a strict preference relation was represented by an open set of s.d.e.u. functions whose extreme points were state-independent, and in extending that relation by imposing $\varepsilon$-tighter weak preferences, it was impossible to retain any agreeing state-dependent utilities that were not convex combinations of agreeing state-independent utilities. A second example shows that this is not always the case under axioms A1–A6: the representing set of state-dependent utilities is not always a convex hull of state-independent utilities, so that nonconstant lotteries may have lower or upper expected utilities that are not supported by any probability/utility pairs. This is not due to the choice of weak vs. strict preference as a primitive, but rather due to an inherent weakness of the usual state-independence axiom in the absence of completeness. A6 serves (only) to extend conditional preferences between constant lotteries to other conditioning events, but, without completeness, there may be few pairs of constant lotteries whose conditional preferences are known, and hence, there may be limited opportunity to apply this axiom.

To illustrate this problem, let there be three states and three consequences, and let $\mathbf{X}$ be the lottery defined by $X(1,0) = X(2,2) = X(3,1) = 1$. That is, $\mathbf{X}$ yields consequences 0, 2 and 1 with certainty in states 1, 2 and 3 respectively. Suppose that all states are judged to have probability at least 0.1, and $\mathbf{X}$ is judged to have an unconditional expected utility of at least 0.5. Furthermore, a coin flip between $\mathbf{X}$ and {consequence 2 if state 1, otherwise $\mathbf{Z}$} is preferred to a coin flip between utility 0.5 and {utility 0.9 if state 1, otherwise $\mathbf{Z}$}, but also a coin flip between $\mathbf{X}$ and {utility 0.1 if state 2, otherwise $\mathbf{Z}$} is preferred to a coin flip between utility 0.5 and {consequence 2 given state 2, otherwise $\mathbf{Z}$}. The common alternative $\mathbf{Z}$ in each comparison is arbitrary by Theorem 1. The basis for $\succsim$ can then be expressed as

(5.1) $$\mathbf{H_E} \succsim \mathbf{H}_{0.1} \quad \text{for } \mathbf{E} = \mathbf{E}_1, \mathbf{E}_2, \mathbf{E}_3,$$

(5.2) $$\mathbf{X} \succsim \mathbf{H}_{0.5},$$

(5.3) $$\tfrac{1}{2}\mathbf{E}_1\mathbf{H}_2 + \tfrac{1}{2}\mathbf{X} \succsim \tfrac{1}{2}\mathbf{E}_1\mathbf{H}_{0.9} + \tfrac{1}{2}\mathbf{H}_{0.5},$$

(5.4) $$\tfrac{1}{2}\mathbf{E}_2\mathbf{H}_{0.1} + \tfrac{1}{2}\mathbf{X} \succsim \tfrac{1}{2}\mathbf{E}_2\mathbf{H}_2 + \tfrac{1}{2}\mathbf{H}_{0.5},$$

where the arbitrary common term $\tfrac{1}{2}(1-\mathbf{E}_1)\mathbf{Z}$ has been suppressed on both sides of (5.3), and likewise $\tfrac{1}{2}(1-\mathbf{E}_2)\mathbf{Z}$ has been suppressed in (5.4). Notice that, by Theorem 1, (5.2) implies $\tfrac{1}{2}\mathbf{E_1Y} + \tfrac{1}{2}\mathbf{X} \succsim \tfrac{1}{2}\mathbf{E}_1\mathbf{Y} + \tfrac{1}{2}\mathbf{H}_{0.5}$ for any $\mathbf{Y}$, from which (5.3) has been constructed by replacing $\mathbf{Y}$ with $\mathbf{H}_2$ on the LHS and replacing it with $\mathbf{H}_{0.9}$ on the RHS. Similarly, (5.2) implies $\tfrac{1}{2}\mathbf{E}_2\mathbf{Y} +$ $\frac{1}{2}\mathbf{X} + \succsim \frac{1}{2}\mathbf{E}_2\mathbf{Y} + \frac{1}{2}\mathbf{H}_{0.5}$, from which (5.4) has been constructed by replacing $\mathbf{Y}$ with $\mathbf{H}_{0.1}$ on the LHS and with $\mathbf{H}_2$ on the RHS. These two new preferences imply that the lower bound on the expected utility of $\mathbf{X}$ among all *probability/utility pairs* agreeing with $\succsim$ must be strictly greater than 0.5, which can be seen as follows. First, an s.d.e.u. function $v$ agrees with (5.3) if and only if $U_v(\mathbf{X}) - U_v(\mathbf{H}_{0.5}) \geq U_v(\mathbf{E}_1\mathbf{H}_{0.9}) - U_v(\mathbf{E}_1\mathbf{H}_2)$, and it agrees with (5.4) if and only if $U_v(\mathbf{X}) - U_v(\mathbf{H}_{0.5}) \geq U_v(\mathbf{E}_2\mathbf{H}_2) - U_v(\mathbf{E}_2\mathbf{H}_{0.1})$. Second, because $U_v(\mathbf{H}_{0.1}) < U_v(\mathbf{H}_{0.9})$, any $v$ must satisfy either $U_v(\mathbf{H}_{0.9}) - U_v(\mathbf{H}_2) > 0$ or $U_v(\mathbf{H}_2) - U_v(\mathbf{H}_{0.1}) > 0$ or both, and if $v$ is a probability/utility pair that agrees with (5.1), it must also satisfy the corresponding conditional inequalities $U_v(\mathbf{E}_1\mathbf{H}_{0.9}) - U_v(\mathbf{E}_1\mathbf{H}_2) > 0$ or $U_v(\mathbf{E}_2\mathbf{H}_2) - U_v(\mathbf{E}_2\mathbf{H}_{0.1}) > 0$ or both, because $\mathbf{E}_1$ and $\mathbf{E}_2$ have positive lower probability. Hence, a probability/utility pair $v$ can agree with (5.1), (5.3) and (5.4) only if $U_v(\mathbf{X}) - U_v(\mathbf{H}_{0.5}) > 0$, that is, $U_v(\mathbf{X}) > 0.5$. In fact, by nonlinear programming, the minimum expected utility of $\mathbf{X}$ among all agreeing probability/utility pairs is achieved at $p(1) = 0.41, p(2) = 0.1, p(3) = 0.49, u(2) = 0.379/0.51 = 0.74314$, and its value is $p(2)u(2) + p(3) = 0.564314$.

However, by direct application of axiom A6 together with A1–A5, the expected utility of $\mathbf{X}$ cannot be determined to be strictly greater than 0.5. To apply A6, it is first necessary to find nonnegative linear combinations of the differences between the LHS's and RHS's of (5.1)–(5.4) that are conditionally constant—that is, of the form $\mathbf{EB}$, where $\mathbf{E}$ is a not-potentially-null event and $\mathbf{B}$ is constant across states. But the search for such conditionally constant terms is constrained here by the presence of a nonconstant term proportional to $\mathbf{X} - \mathbf{H}_{0.5}$ in the differences between LHS's and RHS's of (5.2)–(5.4), as well as by the use of two different conditioning events in (5.3) and (5.4). Furthermore, in order for A6 to "bite," $\mathbf{B}$ needs to have a negative lower expected utility when conditioned on some other event $\mathbf{F}$. The effect of applying A6 will then be to raise this lower expected utility to zero, which shrinks the set of s.d.e.u. functions representing $\succsim$. In the example the few conditionally-constant lottery differences $\mathbf{EB}$ that can be constructed from (5.1)–(5.4) all turn out to satisfy $\mathbf{B} \geq^* 0$, which is uninformative. (This can be determined by solving a sequence of 14 linear programs.) The lower expected utility of $\mathbf{X}$ therefore remains at 0.5 despite the fact that this value is not realized, *or even closely approached*, by any probability/utility pair agreeing with $\succsim$.

This example shows that when preferences are incomplete, axiom A6 is insufficient to guarantee that they are represented by (the convex hull of ) a set of probability/utility pairs. Evidently, a stronger condition is needed, such as follows:

A6* (*Strong state-independence*) If $\mathbf{X}$ and $\mathbf{Y}$ are constant, $\mathbf{X}' \succsim \mathbf{Y}'$, and $\mathbf{H_E} \succsim \mathbf{H}_p$ and $\mathbf{H_F} \precsim \mathbf{H}_q$ with $p > 0$, then

$$\alpha \mathbf{EX} + (1-\alpha)\mathbf{X}' \succsim \alpha \mathbf{EY} + (1-\alpha)\mathbf{Y}'$$

$$\implies \quad \beta\mathbf{FX} + (1-\beta)\mathbf{X}' \succsim \beta\mathbf{FY} + (1-\beta)\mathbf{Y}'$$

for $\beta = 1$ if $\alpha = 1$ and otherwise for all $\beta$ such that $\frac{\beta}{1-\beta} \leq \frac{\alpha}{1-\alpha}\frac{p}{q}$.

In words, A6* requires that whenever a weak preference between two arbitrary lotteries $\mathbf{X}'$ and $\mathbf{Y}'$ is preserved under an objective mixture with two conditional constant lotteries $\mathbf{EX}$ and $\mathbf{EY}$, respectively, and the common conditioning event $\mathbf{E}$ is not potentially null, then the same preference is preserved under a mixture with the same constant lotteries conditioned on any other common event $\mathbf{F}$, where the odds-ratio of the new mixture can be at least as large as the original odds-ratio multiplied by the ratio of the *lower* probability of $\mathbf{E}$ to the *upper* probability of $\mathbf{F}$. In terms of a primal set $\mathcal{B}^{***}$ of preferred directions, A6* requires that if $\mathbf{B}' \in \mathcal{B}^{***}$ and $\mathbf{B}' + \mathbf{EB} \in \mathcal{B}^{***}$, where $\mathbf{B}$ is constant and $\mathbf{E}$ is not potentially null, then $\mathbf{B}' + (p/q)\mathbf{FB} \in \mathcal{B}^{***}$. If $\succsim$ satisfies A1–A6, then it must have an extension that also satisfies A6*, because A1–A6 require the existence of at least one agreeing probability/utility pair, and the additional preferences implied by A6* do not eliminate any agreeing probability/utility pairs. The scale factor $p/q$ ensures precisely that, even in the worst case where $U_v(\mathbf{B}') > 0$, $U_v(\mathbf{B}) < 0$, $p_v(\mathbf{E}) = p$ and $p_v(\mathbf{F}) = q$, $U_v(\mathbf{B}') + p_v(\mathbf{E})U_v(\mathbf{B}) \geq 0 \Longrightarrow U_v(\mathbf{B}') + (p/q)p_v(\mathbf{F})U_v(\mathbf{B}) \geq 0$ for any $v$ that is a probability/utility pair.

The stronger axiom *does* affect the counterexample discussed above. Let $\mathbf{B}' = \mathbf{X} - \mathbf{H}_{0.5}$, in terms of which (5.2), (5.3) and (5.4) are equivalent to $\mathbf{B}' \in \mathcal{B}^{***}$, $\mathbf{B}' + \mathbf{E}_1(\mathbf{H}_2 - \mathbf{H}_{0.9}) \in \mathcal{B}^{***}$ and $\mathbf{B}' + \mathbf{E}_2(\mathbf{H}_{0.1} - \mathbf{H}_2) \in \mathcal{B}^{***}$, respectively. By linear programming, $\mathbf{E}_1$, $\mathbf{E}_2$ and $\mathbf{E}_3$ are found to have identical lower probabilities of 0.1 and upper probabilities of 0.41, 0.4444... and 0.8, respectively. A6* can then be applied to add the following preferred directions to $\mathcal{B}^{***}$: $\mathbf{B}' + (0.1/0.4444...)\mathbf{E}_2(\mathbf{H}_2 - \mathbf{H}_{0.9})$, $\mathbf{B}' + (0.1/0.8)\mathbf{E}_3(\mathbf{H}_2 - \mathbf{H}_{0.9})$, $\mathbf{B}' + (0.1/0.41)\mathbf{E}_1(\mathbf{H}_{0.1} - \mathbf{H}_2)$ and $\mathbf{B}' + (0.1/0.8)\mathbf{E}_3(\mathbf{H}_{0.1} - \mathbf{H}_2)$. These vectors must have nonnegative expected utility under any agreeing s.d.e.u. function, which immediately raises the lower expected utility of $\mathbf{X}$ from 0.5 to 0.54, and further iterations of A6* approach the lower limit of 0.564314 among agreeing probability/utility pairs. However, iteration of A6* via linear programming and vertex enumeration is a cumbersome way to compute probability and utility bounds. The main theorem of this section shows that there is a simpler way:

THEOREM 4. $\succsim$ *satisfies* A1–A5 *and* A6* *if and only if it is represented by a nonempty set* $\mathcal{V}^{***} \in \mathcal{V}^{++}$ *of s.d.e.u. functions that is the convex hull of a set of probability/utility pairs. In particular, if* $\{\mathbf{X}_n \succsim \mathbf{Y}_n\}$ *is a basis for* $\succsim$, *then* $\mathcal{V}^{***}$ *is the convex hull of the set of probability/utility pairs in* $\mathcal{V}^{++}$ *that satisfy* $\{U_v(\mathbf{X}_n) \geq U_v(\mathbf{Y}_n)\}$.

If the basis is finite, the construction of $\mathcal{V}^{***}$ can be carried out as follows. First, form the convex polytope defined by the constraints $\{U_v(\mathbf{X}_n) \geq U_v(\mathbf{Y}_n)\}$, $v \in \mathcal{V}^{++}$. Next, take the intersection of this polytope with the nonconvex set of all probability/utility pairs. (If the latter intersection is empty, the preferences do not satisfy A5 given A1–A4 and A6*: they are incoherent.) Finally, take the convex hull of what remains: this is the set $\mathcal{V}^{***}$.

**6. Discussion.** It has been shown that, in order to obtain a convenient representation of incomplete preferences by sets of probability/utility pairs, it does not suffice merely to delete the completeness axiom from the standard axiomatic framework of Anscombe and Aumann, due to a fundamental weakness of the traditional state-independence axiom in the absence of completeness. The approach here has been to substitute a stronger state-independence postulate (A6*), which still has "bite" in the absence of completeness. This approach follows a common theme in axiomatic rational choice theory, namely, the search for a small number of independent, constructive norms of behavior which collectively turn out to imply the existence of measurable beliefs and values.

Seidenfeld [26] has suggested a modified approach: instead of directly strengthening the state-independence property, simply "compel the missing preferences" by the principle of indirect reasoning, namely, compel the preference $\mathbf{Y} \succ \mathbf{X}$ [$\mathbf{Y} \succsim \mathbf{X}$] whenever the contrary assertion $\mathbf{X} \succsim \mathbf{Y}$ [$\mathbf{X} \succ \mathbf{Y}$] would lead to a violation of the prevailing axioms in light of preferences already known. According to this principle, wherever a weak [strict] preference is "precluded," the opposite strict [weak] preference must be affirmed. In the example of the previous section, it is not directly implied by A1–A6 that $\mathbf{X} \succsim \mathbf{H}_u$ for any $u > 0.5$, yet the same axioms preclude that $\mathbf{H}_u \succsim \mathbf{X}$ for any $u < 0.564314$. If indirect reasoning were invoked, the preferences $\mathbf{X} \succ \mathbf{H}_u$ would immediately be compelled for all $u < 0.564314$. By the same token, Theorem 3 would suffice to fix the lower [upper] expected utility of any lottery at the minimum [maximum] achieved among all agreeing probability/utility pairs, and the stronger state-independence axiom A6* would not be needed for the tighter representation of Theorem 4. (SSK use indirect reasoning in a more limited role to deal with closure of the sets of target utilities in their model, a step which is avoided here through the use of weak rather than strict preference as a primitive.)

Indirect reasoning, otherwise known as *modus tollens* $((p \to q) \wedge \neg q) \to \neg p)$, is a tautological implication of two-valued propositional logic; hence, it applies to preferences that are complete: if either $\mathbf{X} \succsim \mathbf{Y}$ or else $\mathbf{Y} \succ \mathbf{X}$ must be true, then if one is precluded, the other naturally must be affirmed. Indirect reasoning also applies tautologically to incomplete preferences under axioms A1–A5 (only), because according to Theorem 1, the preference $\mathbf{X} \succsim \mathbf{Y}$ merely adds the new direction $\mathbf{X} - \mathbf{Y}$ to the convex cone $\mathcal{B}^*$ of preferred directions, which is to say, it interacts *linearly* with preferences already known. Thus, $\mathbf{X} \succsim \mathbf{Y}$ is precluded precisely when $\mathbf{X} - \mathbf{Y}$ combines linearly with an existing element of $\mathcal{B}^*$ to produce a vector in the open negative orthant, implying that $\mathbf{Y} - \mathbf{X} + \mathbf{Z} \in \mathcal{B}^*$ for some $\mathbf{Z} \in \mathcal{B}^-$, which does entail $\mathbf{Y} \succ \mathbf{X}$. However, the tautology breaks down when the conventional state-independence axiom A6 is added to the mix, because the new direction $\mathbf{X} - \mathbf{Y}$ then interacts *nonlinearly* with an already-known direction of preference if they can be linearly combined to yield a conditionally constant direction that can be propagated to other conditioning events. When the new conditional directions are also added to the cone, a vector in the open negative orthant may result, thus precluding $\mathbf{X} \succsim \mathbf{Y}$ even if it is not already implied that $\mathbf{Y} \succ \mathbf{X}$. Under these conditions, indirect reasoning becomes a logically separate "axiom" of rational choice to which we may or may not wish to subscribe.

The argument in favor of indirect reasoning is that, at the end of the day, the decision maker in the example of Section 5 might face a choice between, say, $\mathbf{X}$ and $\mathbf{H}_{0.54}$. If the *preference* $\mathbf{H}_{0.54} \succsim \mathbf{X}$ is precluded, it would appear as though the *choice* of $\mathbf{H}_{0.54}$ over $\mathbf{X}$ should also be precluded, in which case the implication $\mathbf{X} \succ \mathbf{H}_{0.54}$ would appear to be inescapable as a requirement of rational choice. However, the counterargument is that the point of dropping the completeness assumption is precisely to permit the decision maker to judge some alternatives to be noncomparable, which weakens the link between preference and choice to a one-way implication: a preferred alternative should be chosen where it exists, but a chosen alternative need not be preferred if the alternatives are regarded as noncomparable. In the example, $\mathbf{X}$ and $\mathbf{H}_{0.54}$ are noncomparable under axioms A1–A6 given (only) the preferences (5.1)–(5.4), hence, either may be chosen on that basis, and the choice of $\mathbf{H}_{0.54}$ over $\mathbf{X}$ is in fact supported by some agreeing state-*dependent* utilities that are not ruled out by A6 because of the relatively weak constraints it imposes on the shape of the agreeing set $\mathcal{V}^{**}$. (What is actually precluded is not the preference of $\mathbf{H}_{0.54}$ over $\mathbf{X}$ per se, but rather the additional conditional preferences that would be implied by declaring $\mathbf{H}_{0.54} \succsim \mathbf{X}$ in the presence of A6.) Still, even if indirect reasoning is not an inescapable requirement of rational choice, it can be defended as an axiom that might be useful to enforce when preferences are incomplete, as it is used elsewhere in the foundations of mathematics (e.g., in geometry). However, in this regard, it does not have the elementary character of the other axioms and cannot be tested independently of them, so it is still of interest to ask whether there is some other "direct" postulate that will reduce indirect reasoning once again to a tautology, and that is precisely what A6* accomplishes here.

This paper has considered the case of finite sets of states and consequences. SSK consider the somewhat more general case of a countable set of consequences and a finite partition of a possibly-infinite underlying state space. However, it is often desirable to consider even more abstract settings involving uncountable sets. The case in which the state space is a compact subset of $\Re^n$ and the consequence space is a compact subset of $\Re$ or $\Re^m$ is of particular importance (e.g., in financial economics), where a lottery would naturally be defined by a probability density or cumulative mass at a given consequence in a given state. The axioms of this paper are applicable to such objects, and Theorems 1 and 2 can still be obtained by appropriately generalized separating hyperplane theorems: the primal representation of $\succsim$ is again a convex cone whose directions are differences between more-preferred and less-preferred lotteries, and the dual representation is a convex set of agreeing s.d.e.u. functions; indeed, just such a representation is given by Garcia Del Amo and Ríos Insua [6]. The further restriction to state-independent representations consisting of sets of probability/utility pairs is more problematic but very much of interest in order to reduce the dimensionality of the parameter space in Bayesian robustness applications. Under suitable regularity (e.g., smoothness) conditions, it might be conjectured that a representation of incomplete preferences by sets of state-independent or almost state-independent utilities could still be obtained, since any discretization of the state space and consequence space would have to satisfy Theorems 3 and 4 of this paper. The issue of "missing preferences" would presumably still arise, requiring either an axiom such as A6* of this paper or else the indirect reasoning principle. An analogous set of results in Savage's framework (which requires an infinite state space but does not involve objective probabilities) would appear to be more difficult to obtain, since it does not lend itself to the same techniques of convex analysis that are made possible by objective mixtures, although his state-independence axiom is essentially the same as that of Anscombe–Aumann and might be expected to suffer a similar loss of traction in the absence of completeness.

## APPENDIX

PROOF OF THEOREM 1. Suppose that $\mathbf{X} - \mathbf{Y} \propto \mathbf{X}' - \mathbf{Y}'$, that is, $\alpha(\mathbf{X} - \mathbf{Y}) = (1-\alpha)(\mathbf{X}' - \mathbf{Y}')$ for some $\alpha \in (0,1)$, where $\mathbf{X}, \mathbf{Y}, \mathbf{X}', \mathbf{Y}'$ are lotteries. Then, by the mixture-independence axiom, $\mathbf{X} \succsim \mathbf{Y} \Leftrightarrow \alpha\mathbf{X} + (1-\alpha)\mathbf{X}' \succsim \alpha\mathbf{Y} + (1-\alpha)\mathbf{X}' = \alpha\mathbf{X} + (1-\alpha)\mathbf{Y}' \Leftrightarrow \mathbf{X}' \succsim \mathbf{Y}'$, establishing that there is a collection of vectors $\mathcal{B}^*$ such that $\mathbf{B} \in \mathcal{B}^* \Rightarrow \alpha\mathbf{B} \in \mathcal{B}^* \; \forall \alpha > 0$ and $\mathbf{X} \succsim \mathbf{Y} \Leftrightarrow \mathbf{X} - \mathbf{Y} \in \mathcal{B}^*$. Mixture-independence and transitivity together imply that if $\mathbf{X} \succsim \mathbf{Y}$ and $\mathbf{X}' \succsim \mathbf{Y}'$, then $\alpha\mathbf{X} + (1-\alpha)\mathbf{X}' \succsim \alpha\mathbf{Y} + (1-\alpha)\mathbf{Y}'$, establishing that $\mathcal{B}^*$ is convex. The best/worst axiom establishes that $\mathcal{B}^*$ is nonempty (in particular, it includes $\mathbf{H}_c - \mathbf{H}_0$ and $\mathbf{H}_1 - \mathbf{H}_c$ for $c > 1$), but coherence requires that it not be the whole space (in particular, it excludes $\mathbf{H}_0 - \mathbf{H}_1$), and the continuity axiom implies that it is closed. $\square$

PROOF OF THEOREM 2. Every convex cone $\mathcal{B}^* \subset \mathcal{B}$ is associated with a dual cone $\mathcal{V}^* \subset \mathcal{B}$ such that $\mathbf{B} \in \mathcal{B}^*$ if and only if $U_v(\mathbf{B}) \equiv v^T\mathbf{B} \geq 0$ for all $v \in \mathcal{V}^*$, and coherence requires that $U_v(\mathbf{H}_1 - \mathbf{H}_0) > 0$ for all $v \in \mathcal{V}^*$. Since the elements of $\mathcal{B}^*$ are differences between pairs of lotteries, their components sum to zero within states. Hence, $U_v(\mathbf{B}) \geq 0 \Leftrightarrow U_{v'}(\mathbf{B}) \geq 0$, where $v'(s,c) = \alpha v(s,c) + \beta(s)$ for $\alpha > 0$ and arbitrary $\{\beta(s)\}$, so w.l.o.g. it can be assumed that $v \in \mathcal{V}^+$. $\square$

PROOF OF THEOREM 3. It suffices to show that a preference relation satisfying A1–A6 can be extended to assign a unique state-independent utility to any consequence $c > 1$. This extension can then be performed for $c = 2, 3, \ldots, K$, in turn, to obtain a unique state-independent utility function, for which every agreeing s.d.e.u. function is a probability/utility pair. Therefore, consider some consequence $c > 1$, and let $u(c)$ denote the maximum value of $u$ such that $\mathbf{H}_c \succsim \mathbf{H}_u$—that is, the *greatest lower utility* for consequence $c$. (The maximum exists by virtue of continuity.) Since A1–A5 apply, there exists a cone of preferred directions representing $\succsim$. When A6 also applies, let this cone be denoted by $\mathcal{B}^{**}$: in this case it contains all vectors of the form $\mathbf{E}_s(\mathbf{H}_c - \mathbf{H}_{u(c)})$. Then $\succsim$ can be coherently extended so as to assign $c$ an *upper* utility also equal to $u(c)$. To show this, let $\mathcal{B}_c^{**}$ denote the convex hull of $\mathcal{B}^{**}$ and the rays whose directions are $\mathbf{E}_s(\mathbf{H}_{u(c)} - \mathbf{H}_c)$ for all states $s$, noting that these vectors are opposite in sign to those already known (by A6) to belong to $\mathcal{B}^{**}$. Since $\mathcal{B}_c^{**}$ is also a convex cone, it represents a relation satisfying A1–A4; and since it is formed by adding to $\mathcal{B}^{**}$ a difference between two constant lotteries and all the conditionalizations thereof, the relation it represents also satisfies A6. If $\mathcal{B}_c^{**}$ is also disjoint from the negative orthant $\mathcal{B}^-$, then the relation it represents is coherent and is the desired extension $\succsim_c$. It remains to be shown, then, that $\mathcal{B}_c^{**}$ is disjoint from $\mathcal{B}^-$. Suppose that this is *not* true. Then there must exist $\mathbf{B} \in \mathcal{B}^{**}$ and $\{\beta(s)\} \in [0,1]$ such that

$$\left[-\mathbf{B} - \left(\sum_{s \in S} \beta(s)\mathbf{E}_s(\mathbf{H}_{u(c)} - \mathbf{H}_c)\right)\right]_{\min} > 0$$

$$\iff \left[-\mathbf{B} - \left(\sum_{s \in S}(1-\beta(s))\mathbf{E}_s(\mathbf{H}_c - \mathbf{H}_{u(c)})\right) - (\mathbf{H}_{u(c)} - \mathbf{H}_c)\right]_{\min} > 0$$

$$\iff [-\mathbf{B}' - (\mathbf{H}_{u(c)} - \mathbf{H}_c)]_{\min} > 0,$$

where

$$\mathbf{B}' \equiv \mathbf{B} + \sum_{s \in S}(1-\beta(s))\mathbf{E}_s(\mathbf{H}_c - \mathbf{H}_{u(c)}).$$

Note that $\mathbf{B}' \in \mathcal{B}^{**}$, because $\mathcal{B}^{**}$ includes all vectors of the form $\mathbf{E}(\mathbf{H}_c - \mathbf{H}_{u(c)})$ and their nonnegative linear combinations. Rearrangement of the last inequality yields $(\mathbf{H}_c - \mathbf{H}_{u(c)}) > \mathbf{B}'$, which is equivalent to $(\mathbf{H}_c - \mathbf{H}_{u(c)+\varepsilon}) \geq \mathbf{B}'$ for some $\varepsilon > 0$. Because $\mathbf{H}_c - \mathbf{H}_{u(c)+\varepsilon}$ stochastically dominates $\mathbf{B}' \in \mathcal{B}^{**}$, it follows that $\mathbf{H}_c - \mathbf{H}_{u(c)+\varepsilon} \in \mathcal{B}^{**}$, whence also $\mathbf{H}_c \succsim \mathbf{H}_{u(c)+\varepsilon}$. Thus, $u(c) + \varepsilon$ is a lower utility for $c$, contradicting the assumption that $u(c)$ was the greatest lower utility for $c$ under the relation $\succsim$, therefore $\mathcal{B}_c^{**}$ must be disjoint from $\mathcal{B}^-$. Finally, the same method of extension, namely, adding a preference which asserts that the greatest lower utility is also an upper utility, can be applied to an arbitrary *constant* lottery $\mathbf{X}$ prior to $\{\mathbf{H}_c\}$ to show that the maximum or minimum expected utility of a constant lottery must be achieved by an agreeing probability/utility pair. □

PROOF OF THEOREM 4. The earlier results (based on A1–A6 only) establish that $\succsim$ is represented by a convex cone of preferred directions and dually by a convex set of s.d.e.u. functions containing at least one probability/utility pair. When the stronger A6* holds, let these sets be denoted by $\mathcal{B}^{***}$ and $\mathcal{V}^{***}$, respectively. It must be shown that $\mathcal{V}^{***}$ is the convex hull of a set of probability/utility pairs, which means that, for any lotteries $\mathbf{X}, \mathbf{Y}$, the minimum of $U_v(\mathbf{X}) - U_v(\mathbf{Y})$ over all $v \in \mathcal{V}^{***}$ is achieved at some $v$ which is a probability/utility pair. Suppose that $d = \min_{v \in \mathcal{V}^{***}}(U_v(\mathbf{X}) - U_v(\mathbf{Y}))$. Since $U_v(\mathbf{H}_0) = 0$ and $U_v(\mathbf{H}_u) = u$, it follows that $d$ is the maximum value of $u$ for which

$$(7.1) \qquad U_v(\tfrac{1}{2}\mathbf{X} + \tfrac{1}{2}\mathbf{H}_0) \geq U_v(\tfrac{1}{2}\mathbf{Y} + \tfrac{1}{2}\mathbf{H}_u)$$

for all $v \in \mathcal{V}^{***}$, or, equivalently, the maximum value of $u$ for which $\tfrac{1}{2}\mathbf{X} + \tfrac{1}{2}\mathbf{H}_0 \succsim \tfrac{1}{2}\mathbf{Y} + \tfrac{1}{2}\mathbf{H}_u$, or, equivalently, the maximum value of $u$ for which $\mathbf{X} + \mathbf{H}_0 - \mathbf{Y} - \mathbf{H}_u \in \mathcal{B}^{***}$. To prove that there is a probability/utility pair for which (7.1) holds with equality, it will suffice to show that $\succsim$ has an extension satisfying A1–A6, denoted $\succsim_d$, in which the reverse preference also holds, namely, $\tfrac{1}{2}\mathbf{Y} + \tfrac{1}{2}\mathbf{H}_d \succsim_d \tfrac{1}{2}\mathbf{X} + \tfrac{1}{2}\mathbf{H}_0$. [It is not necessary to enforce the stronger axiom A6* on $\succsim_d$, because if (7.1) holds with equality for all $v$ agreeing with $\succsim_d$, A6 guarantees via Theorem 3 that at least one must be a probability/utility pair.] The primal representation of $\succsim_d$ is constructed as follows. First, define $\mathbf{B}_d = \mathbf{X} + \mathbf{H}_0 - \mathbf{Y} - \mathbf{H}_d$, noting that it is an extreme direction of $\mathcal{B}^{***}$. Then include the reverse preference by adding to the convex cone $\mathcal{B}^{***}$ the ray whose direction is $-\mathbf{B}_d$. Then apply A1–A3, which convexifies the set, yielding an expanded cone. (A4 and A5 are automatically satisfied thus far.) Finally, apply A6 (the special case of A6* with $\alpha = 1$), which means to (i) determine which new conditional directions of the form $\mathbf{F}_i \mathbf{B}_i$, where $\mathbf{F}_i$ is a not-potentially-null event and $\mathbf{B}_i$ is constant, are in the expanded cone; (ii) then propagate the same conditional preferences to all possible conditioning events by adding the directions $\{\mathbf{E}_s \mathbf{B}_i\}$ for all $s$ and $i$; and (iii) take the convex hull again. Thus, the primal cone consists of all vectors of the form $\mathbf{B}_0 - \gamma_d \mathbf{B}_d + \sum_{s,i} \gamma_{si} \mathbf{E}_s \mathbf{B}_i$, where $\mathbf{B}_0 \in \mathcal{B}^{***}$, $\gamma_d \geq 0$ and $\gamma_{si} \geq 0$. Suppose that one such vector is in the open negative orthant $\mathcal{B}^-$, in violation of A5, that is, suppose that $\mathbf{B}_0 - \gamma_d \mathbf{B}_d + \sum_{s,i} \gamma_{si} \mathbf{E}_s \mathbf{B}_i = \mathbf{Z}$, or, equivalently,

$$(7.2) \qquad \mathbf{B}_0 + \sum_{s,i} \gamma_{si} \mathbf{E}_s \mathbf{B}_i = \gamma_d \mathbf{B}_d + \mathbf{Z},$$

where $\mathbf{B}_0 \in \mathcal{B}^{***}$ and $\mathbf{Z} \in \mathcal{B}^-$. It will be shown that this leads to a contradiction.

Because the conditional direction $\mathbf{F}_i \mathbf{B}_i$ is in the convex hull of $\mathcal{B}^{***}$ *plus* the additional direction $-\mathbf{B}_d$, but not in $\mathcal{B}^{***}$ alone, it satisfies $\mathbf{B}_d + \mathbf{F}_i \mathbf{B}_i \in \mathcal{B}^{***}$. Because $\mathbf{B}_d \in \mathcal{B}^{***}$ and each $\mathbf{F}_i$ is not-potentially-null and the original relation $\succsim$ satisfies A6*, there exist positive constants $\{\delta_{si}\}$ such that $\mathbf{B}_d + \delta_{si} \mathbf{E}_s \mathbf{B}_i \in \mathcal{B}^{***}$ for all $i$ and $s$. Multiplying through by $\gamma_{si}/\delta_{si}$, summing and invoking the convexity of $\mathcal{B}^{***}$ yields $(\sum_{s,i} \gamma_{si}/\delta_{si}) \mathbf{B}_d + \sum_{s,i} \gamma_{si} \mathbf{E}_s \mathbf{B}_i \in \mathcal{B}^{***}$, whence also $\mathbf{B}_0 + (\sum_{s,i} \gamma_{si}/\delta_{si}) \mathbf{B}_d + \sum_{s,i} \gamma_{si} \mathbf{E}_s \mathbf{B}_i \in \mathcal{B}^{***}$. Adding $(\sum_{s,i} \gamma_{si}/\delta_{si}) \mathbf{B}_d$ to both sides of (7.2) and comparing, it follows that $(\gamma_d + \sum_{s,i} \gamma_{si}/\delta_{si}) \mathbf{B}_d + \mathbf{Z} \in \mathcal{B}^{***}$. Multiplication by $k = (\gamma_d + \sum_{s,i} \gamma_{si}/\delta_{si})^{-1}$ yields $\mathbf{B}_d + k\mathbf{Z} \in \mathcal{B}^{***}$, where $k\mathbf{Z} \in \mathcal{B}^-$, contradicting the assumption that $\mathbf{B}_d$ was an extreme direction of $\mathcal{B}^{***}$ and constructively proving the existence of $u > d$ for which $\frac{1}{2}\mathbf{X} + \frac{1}{2}\mathbf{H}_0 \succsim \frac{1}{2}\mathbf{Y} + \frac{1}{2}\mathbf{H}_u$. $\square$

**Acknowledgments.** The author is grateful to David Ríos Insua, Teddy Seidenfeld, James E. Smith, Peter Wakker and two anonymous referees for comments on earlier drafts.

THE FUQUA SCHOOL OF BUSINESS
DUKE UNIVERSITY
DURHAM, NORTH CAROLINA 27708-0120
USA
E-MAIL: robert.nau@duke.edu
URL: www.duke.edu/~rnau